\documentclass[preprint,11pt]{amsart}
\usepackage{amscd,amssymb,amsopn,amsmath,amsthm,graphics,amsfonts,enumerate,verbatim,calc}
\usepackage[dvips]{graphicx}
\usepackage[colorlinks=true,linkcolor=red,citecolor=blue]{hyperref}
\newcommand{\ra}{\rightarrow}

\newcommand{\expxy}{\exp_{x}^{-1}{y}}
\newcommand{\expyx}{\exp_{y}^{-1}{x}}
\newcommand{\norm}{\|}

\newcommand{\supp}{{\rm{Supp\; }}}

\newcommand{\hess}{{\rm{Hess }}}

\newcommand{\se}{{\rm{sec}}}
\newcommand{\im}{{\rm{Im}}}
\newcommand{\dom}{{\rm{dom}}}

\newcommand{\li}{\langle}%%$left inner product
\newcommand{\ri}{\rangle}%%%right inner product

\newcommand{\RR}{\mathbb{R}}

\newtheorem{theorem}{Theorem}[section]
\newtheorem{corollary}[theorem]{Corollary}
\newtheorem{lemma}[theorem]{Lemma}

\newtheorem{proposition}[theorem]{Proposition}

\newtheorem{definition}[theorem]{Definition}

\theoremstyle{plain}

\newcommand{\blu}[1]{{\color{black}{#1}}}

\numberwithin{equation}{section}
\begin{document}
\title[Weak geodesics on prox-regular sets]{Weak geodesics on prox-regular subsets of Riemannian manifolds}
\author[J. Ferrera, M. R. Pouryayevali, and H. Radmanesh]{Juan Ferrera, Mohamad R. Pouryayevali, and Hajar Radmanesh}
\address{IMI, Departamento de An\'{a}lisis Matem\'{a}tico, Facultad Ciencias Matem\'{a}ticas, Universidad Complutense, 28040, Madrid, Spain}
\email{ferrera@mat.ucm.es}
\address{Department of Pure Mathematics, Faculty of Mathematics and Statistics, University of Isfahan, Isfahan,
81746-73441, Iran}
\email{pourya@math.ui.ac.ir} \email{h.radmanesh@sci.ui.ac.ir}

\maketitle

%%%%%%%%%%%%%%%%%%%%%%%%%%%%%%%%%%%%%%%%%%%%%%%%%%%%%%%%%%%%%%%%%%%%%%%%%%%%%%%%%%%%%%%%%%%%%%%%%%%%%%%%%%%%%%%
%%%%%%%%%%%%%%%%%%%%%%%%%%%%%%%%%%%%%%%%%%%%%%%%%%%%%%%%%%%%%%%%%%%%%%%%%%%%%%%%%%%%%%%%%%%%%%%%%%%%%%%%%%%%%%%
\begin{abstract}
We give a definition of weak geodesics  on  prox-regular  subsets of  Riemannian manifolds as continuous  curves with some weak regularities. Then obtaining a suitable Lipschitz constant of the projection map,  we characterize weak geodesics on a prox-regular set  with assigned end points as viscosity critical points of the energy functional.
\end{abstract}

%%%%%%%%%%%%%%%%%%%%%%%%%%%%%%%%%%%%%%%%%%%%%%%%%%%%%%%%%%%%%%%%%%%%%%%%%%%%%%%%%%%%%%%%%%%%%%%%%%%%%%%%%%%%%%%
%%%%%%%%%%%%%%%%%%%%%%%%%%%%%%%%%%%%%%%%%%%%%%%%%%%%%%%%%%%%%%%%%%%%%%%%%%%%%%%%%%%%%%%%%%%%%%%%%%%%%%%%%%%%%%%
\section{Introduction}\label{sec1}
A classical topic in differential geometry and global nonlinear analysis is the study of geodesics on a
Riemannian manifold $M$ without boundary. Considering the paths with assigned extreme points, the problem of
existence and multiplicity of geodesics on manifolds without boundary was studied in the classical
works \cite{Schwartz,Serre}. Indeed, it was proved that they are critical points of the energy functional on the smooth manifold $X$ of the admissible paths and by means of Morse and  Lusternik-Schnirelman  theory, the multiplicity results are obtained.

In the case in which the manifold $M$ has boundary, even if $M$ is smooth, different kind of irregularities may be developed. For example, the natural domain of the energy functional, i.e., the Sobolev space of $W^{1,2}$-paths on $M$ has no more the structure of a Hilbert manifold and there is no uniqueness for the Cauchy problem. Moreover, in this case the geodesics are not in general $C^2$, however they are differentiable curves with locally Lipschitz derivative. To obtain the results regarding this issue, the theory
of critical points and gradient flows  for some lower semicontinuous functions \blu{is employed}; see \cite{Bishop,Marino,Scolozzi,Wolter}.

The above manifolds are basically smooth and it is natural to study manifolds with a certain degree of irregularity.
Various extensions have been considered in this regards, for instance one can refer to conical manifolds.
A conical manifold $M$ is a complete $m$-dimensional $C^0$ submanifold of $\mathbb{R}^n$ which is everywhere smooth, except for a finite set of points, see
\cite{Ghimenti}. Another development in which the manifold $M$ is the closure of a bounded open subset of $\mathbb{R}^n$ with Lipschitz boundary has been started in  \cite{Morbini}. Then  for this case, a new definition of geodesic on a general subset $M$ of $\mathbb{R}^n$ which is related to the nonsmooth critical point theory \blu{developed in \cite{Marzocchi}, is used. The intrinsic case in which $M$ is the
closure of a bounded open subset of a differentiable manifold $N$ has been studied more recently in \cite{Lancelotti}.}

Another natural development was provided in \cite{Canino,Canino1,Canino2}, where geodesics on
certain nonsmooth sets of $\mathbb{R}^n$, called $p$-convex (or $\varphi$-convex) sets,
are considered. In spite of the lack of regularity in the
set $M$, using a new definition of geodesics in the framework of Sobolev spaces the author characterized these geodesics as critical points of the energy functional on a suitable path space. The class of $p$-convex sets includes submanifolds (possibly with boundary) of class $C^{1,1}_{loc}$, images under $C^{1,1}_{loc}$-diffeomorphism of convex sets, but also subsets which are not topological manifolds, although they
are absolute neighborhood retracts. In particular, it contains subsets with corners of convex type and concave parts of class $C^2$. In \cite{Morbini}, the authors  proved that their notion of geodesic agrees with that of \cite{Canino2}, when $M$ is a $C^2$-submanifold of $\mathbb{R}^n$, possibly with boundary.

The notions of $\varphi$-convexity (as titled $p$-convexity)
 and prox-regularity of sets were introduced in \cite{Colombo1} and \cite{Poliquin}, respectively.
In \cite{Barani} the concept of $\varphi$-convex sets was extended  to Hadamard manifolds and it was shown that if $S$ is a $\varphi$-convex subset of an infinite-dimensional Hadamard manifold $M$, then there exists a neighborhood $U$ of $S$ in $M$ such that the metric projection $P_S: U\to S$ is single-valued and locally Lipschitz.
On the other hand, in \cite{Hosseini}  the notion of prox-regular sets  was introduced on Riemannian manifolds as a subclass of regular sets.
Prox-regular sets have significant applications in the theory of Moreau sweeping process, crowd motion and  second order analysis; see, for instance, \cite{Maury,Tanwani}.

In \cite{convex} we proved that   the two classes of $\varphi$-convex sets and prox-regular sets coincide in the setting of finite-dimensional Riemannian manifolds. Moreover, in \cite{Minimizing} we verified that for a prox-regular subset $S$ of a Riemannian manifold $M$, $P_S$ is a locally Lipschitz retraction from a neighborhood of $S$.

In \cite{Minimizing},  the subject of minimizing curves on a prox-regular set $S\subset M$  with $C^2$ boundary was considered. In this paper,  we employed an adapted variational technique and by applying the first variation formula, we obtained a necessary condition for an admissible  curve to be minimizing on $S$. Indeed, this curve is a piecewise $C^2$ curve and it has the property that
\[
D_t\dot{\gamma}(t)\in N^P_S\left(\gamma(t)\right),
 \]
for every $t\in [a,b]$ except for finitely many points, where $N^P_S(x)$ denotes the proximal normal cone at $x\in S$. When the prox-regular set $S$ does not possess a $C^2$ boundary, the problem becomes more complicated and the first variation formula can no longer be applied. Therefore motivated by \cite{Canino2}, we want to study weak geodesics on a prox-regular set $S$ in the intrinsic case where $S$ is a  subset of a Riemannian manifold $M$. Since a prox-regular set $S$ is not necessarily smooth, a geodesic on $S$ needs to be define with some weaker regularities as a curve in the Sobolev space $W^{2,2}(I,M)$. The aim of this paper is to demonstrate that the so-called weak geodesics on $S$ are critical points of the energy functional. In Section \ref{sec2} we review some preliminary
concepts and results from nonsmooth analysis and Riemannian manifolds. In Section \ref{sec3} we give a definition of weak geodesics on $S$ with assigned end points and also we define a constraint  minimization problem using the energy functional. Moreover,  we obtain a suitable Lipschitz constant of the projection map and we prove some auxiliary theorems. Section \ref{sec4} is devoted to the study of critical points of the energy functional and  \blu{the main results of the paper are proved}

%%%%%%%%%%%%%%%%%%%%%%%%%%%%%%%%%%%%%%%%%%%%%%%%%%%%%%%%%%%%%%%%%%%%%%%%%%%%%%%%%%%%%%%%%%%%%%%%%%%%%%%%%%%%%%%%%%
%%%%%%%%%%%%%%%%%%%%%%%%%%%%%%%%%%%%%%%%%%%%%%%%%%%%%%%%%%%%%%%%%%%%%%%%%%%%%%%%%%%%%%%%%%%%%%%%%%%%%%%%%%%%%%%%%%
\section{Preliminaries and notations}\label{sec2}
Let us recall some notions of Riemannian manifolds and nonsmooth analysis; see, e.g.,
\cite{Azagra1,Azagra2,Clark book a,Do Carmo,sakai}. Throughout this paper, $I$ denotes the closed interval $[0,1]$ and $(M,g)$ is an $n$-dimensional Riemannian manifold endowed with a Riemannian metric $g_x=\li .,.\ri_x$ on each tangent space $T_xM$ and  $\nabla$ is the Riemannian connection of $g$.   For every $x, y \in M$, the
Riemannian distance from $x$ to $y$ is denoted by $d(x,y)$. Moreover, $B(x,r)$ and $\overline{B}(x,r)$ signify the
open and closed metric ball centered at $x$ with radius $r$, respectively. For a smooth curve
$\gamma : I \rightarrow M$ and $t_0,t\in I$, the notation  $L^{\gamma}_{t_0,t}$ is used for the parallel
transport along $\gamma$ from $\gamma(t_0)$ to $\gamma(t)$. When $\gamma$ is the unique minimizing geodesic
joining $y$ to $x$, we use the notation $L_{y,x}$.
%Furthermore for a smooth vector field $X$ along $\gamma$, $D_t X$ is the covariant derivative of $X$ along $\gamma$.

For $x\in M$, let $r(x)$ be the convexity radius at $x$, then the function $x\mapsto r(x)$ from
$M$ to $\mathbb{R}^{+} \cup \{+\infty \}$ is continuous; see \cite{sakai}.  The map $\exp_x : U_x \to M$ will stand for the exponential map
at $x$, where $U_x $ is an open subset of the tangent space $T_xM$
containing $0_{x} \in T_xM$.
 Note that if $x$ and $y$ belong to a convex set, then both
$\expxy$ and $\expyx$ are defined and
\[
  \norm \expxy \norm = d(x,y) = \norm \expyx \norm, \qquad L_{y,x}\left(\exp_y^{-1}x\right) = - \exp_x^{-1}y.
\]
Moreover, for a fixed point $z\in M$, the function $\psi:M\ra \RR$ defined by $\psi(x)=d^2(x,z)$ is $C^{\infty}$ on
any convex neighborhood  of $z$ and for every $x$ in a convex neighborhood of $z$,  $\nabla \psi(x)= - 2\exp_x^{-1}z$.

For $x\in M$, let $B\left(x,R\right)$ be a convex ball with compact closure and $\delta,\Delta$ be such that $\delta\leq\se\leq \Delta$ for all sectional curvatures of $M$ on $B\left(x,R\right)$. Then using Rauch's theorem, it can be derived that for any $0<r\leq R$, $\exp_x^{-1}$ is Lipschitz on $B(x,r)$ with the Lipschitz constant $k=k(r)$ defined by
\[
 k=\left\{\begin{array}{lr}
           \frac{r\sqrt{\Delta}}{\sin\left(r\sqrt{\Delta}\right)} & \Delta > 0\\
           1 & \Delta\leq 0
           \end{array}\right..
\]

%%%%%%%%%%%%%%%%%%%%%%%%%%%%%%%%%%%%%%%%%%%%%%%%%%%%%%%%%%%%%%%%%%%%%%%
Let $S$ be a nonempty closed subset of  $M$. Recall that the distance function to $S$ is  $d_S(z)= \inf_{x\in S} \, d(x,z)$ and the metric projection to $S$, denoted by $P_S$, is defined by
\[
P_S(z)= \left\{ x\in S : d_S(z)=d(x,z)\right \}\quad \forall z\in M.
\]
The proximal normal cone to $S$ at $x\in S$, is denoted by $N^{P}_{S}(x)$ and $\xi\in N^{P}_{S}(x)$ if and only if there exists  $\sigma>0$ such that
\[
\li \xi, \expxy \ri \leq\sigma\;d^2(x,y),
\]
for every $y\in U\cap S$, where $U$ is a convex neighborhood of $x$. Moreover, the Bouligand  tangent cone to $S$ at $x$ is defined as
\[
T^B_S(x):=\left\{\lim_{i\ra \infty}\frac{\exp_x^{-1} z_i}{t_i} : z_i\in U\cap S, z_i\ra x\  \textrm{and}\
t_i\downarrow 0\right\},
\]
where $U$ is a convex neighborhood of $x$ in $M$.

Let $f:M\ra (-\infty,+\infty]$ be a lower semicontinuous function and $x\in \dom(f):=\left\{y\in M:f(y)<\infty\right\}$. The viscosity (or Fr\'{e}chet) subdifferential of $f$ at $x$, denoted by $D^-f(x)$, is the set
\[
D^-f(x):=\left\{dg(x): g\in C^1(M,\mathbb{R}), f-g \ \text{attains a local minimum at}\ x\right\}.
\]
Using \cite[Theorem 4.3]{Azagra2}, $\xi \in D^-f(x)\subset T_xM$ if and only if
\[
\liminf_{v\ra 0}\frac{\left(f\circ \exp_x\right)(v)-f(x)-\li \xi,v\ri_x}{\norm v\norm}\geq 0.
\]
It is worth mentioning that if $f$ has a local minimum at $x$, then $0\in D^-f(x)$.
%%%%%%%%%%%%%%%%%%%%%%%%%%%%%%%%%%%%%%%%%%%%%%%%%%%%%%%%%%%%%%%%%%%%%%%%%%%%%

Let us now take a brief look at the subject of prox-regular sets and present some of their properties. A closed subset $S$ of $M$ is prox-regular at $\bar{x}\in S$ if there exist $\varepsilon>0$ and $\sigma>0$ such that $B(\bar {x}, \varepsilon)$ is convex and for every  $x\in S\cap B(\bar{x}, \varepsilon)$ and $v\in N^P_S(x)$ with $\norm v\norm <\epsilon$,
  \[
    \li v, \exp_x^{-1}y \ri \leq\sigma\;d^2(x,y) \quad \forall\  y\in S\cap B(\bar{x}, \epsilon).
  \]
Moreover, $S$ is called prox-regular if it is prox-regular at each point of $S$; for more details, see \cite{Hosseini}.

In \cite[Theorem 3.4]{convex}, it was shown that for every prox-regular subset $S$ of $M$ there exists a continuous function $\varphi : S\rightarrow [0,\infty)$ such that $S$ is  $\varphi$-convex.
Recall that a nonempty closed subset $S\subset M$ is called $\varphi$-convex if for every $x\in S$ and $v\in N^P_S(x)$
  \[
   \li v,\exp_x^{-1}y \ri \leq \varphi(x)\norm v\norm d^2(x,y),
   \]
for every $y\in U\cap S$, where $U$ is a convex neighborhood of $x$. Since we need to utilize the function $\varphi$, we prefer to work with $\varphi$-convex sets.

In \cite{Minimizing}, we proved that for a closed $\varphi$-convex set $S$, the metric projection $P_S$ is locally Lipschitz on an open set containing  $S$. Moreover,  $P_S$ is directionally differentiable at each point $x\in S$ and for every $x\in S$ and $v\in T_{x}M$ we have
\begin{equation}\label{direc}
  \lim_{t\ra 0^+}\frac{\exp_{x}^{-1}\left(P_S\left(\exp_{x}(tv)\right)\right)}{t}=P_{T^B_S\left(x\right)}(v),
\end{equation}
where $P_{T^B_S\left(x\right)}$ denotes the metric projection to the Bouligand tangent cone $T^B_S\left(x\right)$.

%%%%%%%%%%%%%%%%%%%%%%%%%%%%%%%%%%%%%%%%%%%%%%%%%%%%%%%%%%%%%%%%%%%%%%%%%%%%%%%%%%%%%%%%%%%%%%%%%%%%%%%%%%%%%%%%%
%%%%%%%%%%%%%%%%%%%%%%%%%%%%%%%%%%%%%%%%%%%%%%%%%%%%%%%%%%%%%%%%%%%%%%%%%%%%%%%%%%%%%%%%%%%%%%%%%%%%%%%%%%%%%%%%%%%

%\section{Topological properties of the energy functional}\label{sec3}
\section{Lipschitz constant of projection map}\label{sec3}
Our first task in this section is to define weak geodesics on a prox-regular set and the energy functional on a suitable
constraint and to study the lower semicontinuity of the energy functional.

In the case in which the prox-regular set $S$ has  $C^2$ boundary, we were able to obtain a necessary condition for a curve $\gamma$ to be a minimizing curve between its endpoints in $S$, see \cite[Theorem 6]{Minimizing}. In this situation, an admissible curve is a piecewise $C^2$ curve $\gamma:I\ra M$ with nonzero derivatives  that is entirely in $S$. When $S$ is an arbitrary prox-regular set without any smoothness assumption on its boundary, the set of admissible curves needs to be considered more broadly. Since a  prox-regular set is locally Lipschitz path connected, an admissible curve can be chosen among curves belonging to the Sobolev space $W^{1,2}(I,M)$ or so-called $H^1$-curves. An $H^1$-curve $\gamma:I\ra M$ can be considered as an absolutely continuous curve for which $\dot{\gamma}(t)$ exists for almost all $t\in I$ and $\int_{0}^{1}\norm \dot{\gamma}(t)\norm^2 dt<\infty$, see \cite{Klingenberg}.

In order to remind the Sobolev space of  manifold valued curves, we present the following proposition which is a slight modification of \cite[Lemma B.5]{Wehrheim}. By $C^0(I,M)$ we denote the space of continuous curves, endowed with the metric $d_{\infty}\left(\gamma,\eta\right)=\sup_{t\in I}d\left(\gamma(t),\eta(t)\right)$ and $C^{\infty}(I,M)$ denotes the set of smooth curves. For more details about Sobolev spaces, see for instance \cite{Adams}.

Let $\left(U_{\alpha}, \phi_{\alpha}\right)_{\alpha\in A}$ be an atlas on $M$ and $\Phi: M\ra \mathbb{R}^{2n+1}$ be an embedding which exists by the Whitney theorem.
\begin{proposition}\label{prop1}
  Let $k\in \mathbb{N}$ and $1\leq p<\infty$ be such that $kp>1$. Then for $\gamma\in C^0(I,M)$ the following statements are equivalent:
  \begin{itemize}
    \item[(i)] $\phi_{\alpha}\circ\gamma\in W^{k,p}\left(\gamma^{-1}\left(U_{\alpha}\right), \mathbb{R}^{n}\right)$ for all $\alpha$;
    \item[(ii)] $\Phi\circ \gamma\in W^{k,p}\left(I,\mathbb{R}^{2n+1}\right)$;
    \item[(iii)] $\gamma=\exp_c(V)$ for some $c\in C^{\infty}(I,M)$ and $V\in W^{k,p}\left(I, \gamma^{-1}TM\right)$,
  \end{itemize}
  where $\gamma^{-1}TM$ denotes the pullback bundle of $TM$ by $\gamma$.
\end{proposition}

According to the previous proposition, the Sobolev space $W^{k,p}\left(I,M\right)$ of curves on $M$  is defined as the set of continuous curves that satisfy these equivalent statements. In the case $k=0$ and $p=2$, we consider the space $L^2(I,M)$ as
\[
L^2(I,M):=\left\{\gamma:I\ra M : \Phi\circ \gamma\in L^2\left(I,\mathbb{R}^{2n+1}\right)\right\},
\]
with the topology given by the following convergence criteria:
\[
\gamma_i\ra \gamma \ \text{in}\ L^2(I,M)\quad\Leftrightarrow\quad \Phi\circ\gamma_i\ra \Phi\circ\gamma \ \text{in}\ L^2\left(I,\mathbb{R}^{2n+1}\right).
\]

%Recently, the Sobolev space of maps between manifolds has been widely studied in different contexts, see for instance \cite{Michor,Convent1,Convent2,Urakawa,Wehrheim}.

The Sobolev space $W^{1,2}\left(I,M\right)$ usually denotes by $H^1(I,M)$ and as known it is a Hilbert manifold. Moreover, the Sobolev space $H^1\left(I, \gamma^{-1}TM\right)$ is a Hilbert space containing all continuous vector fields $V:I\ra M$ along $\gamma$ with the properties that
\[
\int_{0}^{1}\norm V(t)\norm^2dt<\infty\quad\text{and}\quad \int_{0}^{1}\norm D_tV(t)\norm^2dt<\infty.
\]
Note that for $V\in H^1\left(I, \gamma^{-1}TM\right)$, the weak covariant derivative of $V$, denoted by $D_t V$, is defined as
\[
D_t V:=k_{TM}\circ DV,
\]
where $DV:I\ra \mathbb{R}\otimes TTM$ is the weak derivative of $V$ and $k_{TM}:TTM\ra TM$ is the connection map on the tangent bundle $TM$. Since $\norm D_t V(t)\norm_{\gamma(t)}\leq \norm DV(t)\norm_{\gamma(t)}$, we have $D_t V\in L^2\left(I, \gamma^{-1}TM\right)$; see \cite{Convent1,Convent2}.

Now let $S$ be a  closed and connected  prox-regular subset of $M$ and $U$ be an open neighborhood of $S$ on which $P_S$ is single-valued and locally Lipschitz.  In order to obtain a necessary condition for a curve $\gamma$ to be a minimizing curve between its endpoints, we consider the set of admissible curves as follows:
\[
\mathcal{A}=\mathcal{A}_{x,y}:=\left\{\gamma\in H^1(I,M) : \gamma(t)\in S, \ \forall t\in I, \gamma(0)=x, \gamma(1)=y\right\}.
\]
Since $S$ is locally Lipschitz path connected, $\mathcal{A}$ is nonempty. Note that without considering the constraint ``$\gamma(t)\in S, \ \forall t\in I$'', the set $\mathcal{A}_{x,y}$ is a submanifold of $H^1(I,M)$ and with the distance deduced from the Riemannian metric on $H^1(I,M)$, $\mathcal{A}_{x,y}$ is a complete metric space, see \cite{Klingenberg}.

Motivated by \cite[Theorem 6]{Minimizing} and \cite{Canino2}, we introduce  a
weak geodesic on $S$ as follows:
\begin{definition}
   A continuous curve $\gamma : I\ra M$,  is called a weak geodesic on $S$ joining $x,y\in S$ if $\gamma(0)=x$, $\gamma(1)=y$ and
   \begin{itemize}
    \item [(a)] $\gamma(t)\in S \quad \forall t\in I$;
    \item [(b)] $\gamma\in W^{2,2}(I,M)$;
    \item [(c)] $D_t\dot{\gamma}(t)\in N^P_S\left(\gamma(t)\right) \quad a.\ e.\ t\in I$.
   \end{itemize}
\end{definition}
It is worth mentioning that when $S$ is considered to be all of $M$ or $S$ is a $C^2$-submanifold of $M$ possibly with boundary, this definition corresponds to the usual definition of geodesics.

 We intend  to characterize weak geodesics on $S$ as nonsmooth critical  points  of the energy functional defined on  the space $L^2(I,M)$. We define the energy functional as follows,
\[
\hspace{-3.8cm} f: L^2(I,M)\ra \mathbb{R}\cup \{+\infty\}
\]

\[
f(\gamma)=\left\{
\begin{array}{lr}
  \frac{1}{2} \int_0^1 \norm \dot{\gamma}(t)\norm ^2 dt & \gamma\in \mathcal{A} \\
  & \\
  +\infty & \gamma\in L^2(I,M)\setminus\mathcal{A},
\end{array}\right.
\]
and we consider the optimization problem $\min_{\gamma\in L^2(I,M)}f$.
%It is known that $\mathcal{A}$ has not a natural structure of manifold and $f$ is not regular. These facts suggest that the more natural way to deal with this problem is to use of nonsmooth techniques on Riemannian manifolds.
In the following proposition, we present some topological properties of  $\mathcal{A}$  and the energy functional $f$.

\begin{proposition}
  The admissible set $\mathcal{A}$ is closed in $L^2(I,M)$ and the energy functional $f$ is lower semicontinuous.
\end{proposition}
\begin{proof}
  In order to use the continuity of $P_S$, we consider a subset $\mathcal{D}$ of $H^1(I,M)$ as
  \[
  \mathcal{D}=\left\{\gamma\in H^1(I,M) : \gamma(t)\in U, \ \forall t\in I\right\}.
  \]
  Using (iii) of Proposition \ref{prop1}, it can be proved that $\mathcal{D}$ is open in $H^1(I,M)$. On the other hand, since  $\norm \Phi\circ\gamma\norm_{L^2}\leq \norm \Phi\circ\gamma\norm_{H^1}$ for every $\gamma\in H^1(I,M)$, it follows that $\mathcal{D}$ is open in $L^2(I,M)$.

  We now define the functional map $P_S:\mathcal{D}\ra \mathcal{D}$ by $P_S\left(\gamma\right)(t):=P_S\left(\gamma(t)\right)$ for all $t\in I$. Since the metric projection $P_S$ is locally Lipschitz on $U$, the map $P_S:\mathcal{D}\ra \mathcal{D}$ is well defined. We claim that $P_S:\mathcal{D}\ra \mathcal{D}$ is continuous. Indeed, let $\gamma_n$ be a sequence in $\mathcal{D}$ such that  $d_{\infty}\left(\gamma_n,\gamma\right)\ra 0$. Since $P_S$ is locally Lipschitz on $U$ and $\im\left(\gamma\right)$ is compact, there exists $l>0$ such that
  \[
  d\left(P_S\left(\gamma_n(t)\right),P_S\left(\gamma(t)\right)\right)\leq ld\left(\gamma_n(t),\gamma(t)\right)\leq ld_{\infty}\left(\gamma_n,\gamma\right),
  \]
  for all $t\in I$. This implies that $d_{\infty}\left(P_S\left(\gamma_n\right),P_S\left(\gamma\right)\right)\ra 0$ and then $P_S\left(\gamma_n\right)\ra P_S\left(\gamma\right)$ in $L^2(I,M)$.

  We now define the map $g:\mathcal{D}\ra \mathbb{R}\times M\times M$ as
  \[
  g(\gamma):=\left(d_{\infty}\left(\gamma,P_S\left(\gamma\right)\right),\gamma(0),\gamma(1)\right),
  \]
  hence $g$ is continuous and then $\mathcal{A}=g^{-1}\{(0,x,y)\}$ is closed.
  Therefore the characteristic function $1_{\mathcal{A}}$ of $\mathcal{A}$  is lower semicontinuous. Then as a product of a nonnegative continuous function and a lower semicontinuous function, $f$ is lower semicontinuous, because $f=E\times 1_{\mathcal{A}}$ where $E:H^1(I,M)\ra \mathbb{R}$ is a smooth functional (see \cite{Klingenberg}) defined by
  \[
  E\left(\gamma\right):=\frac{1}{2}\int_{0}^{1}\norm \dot{\gamma}(t)\norm^2dt.
  \]
\end{proof}

%%%%%%%%%%%%%%%%%%%%%%%%%%%%%%%%%%%%%%%%%%%%%%%%%%%%%%%%%%%%%%%%%%%%%%%%%%%%%%%%%%%%%%%%%%%%%%%%%%%%%%%%%%%%%%%%%%
%%%%%%%%%%%%%%%%%%%%%%%%%%%%%%%%%%%%%%%%%%%%%%%%%%%%%%%%%%%%%%%%%%%%%%%%%%%%%%%%%%%%%%%%%%%%%%%%%%%%%%%%%%%%%%%%%%
%\section{Lipschitz constant of projection map}\label{sec4}
\blu{We now improve} the Lipschitz constant of $P_S$ obtained in \cite{Minimizing} and then we prove some auxiliary theorems that are used in the proof of main results of the paper in Section \ref{sec4}.

Let $S$ be a nonempty, closed and $\varphi$-convex subset of $M$ and using \cite[Theorm 4]{Minimizing}, let $U$ be an open set containing $S$ such that $P_S$ is single-valued and locally Lipschitz on $U$.
For simplicity, we use the notation $\phi_x:=\exp_x^{-1}$.
\begin{lemma}\label{Lemma1}
  For every $x_0\in M$ there exist $R > 0$ and $C > 0$ such that the Lipschitz
  constant $A$ of
  \[
   L_{y,x}\circ\phi_{y} - \phi_{x} : B\left(x_0, R\right) \ra T_xM,
  \]
  satisfies $A \leq Cd(x, y)$ for every $x, y \in B\left(x_0, R\right)$.
\end{lemma}
\begin{proof}
 Let $x_0\in M$ and $\mathcal{U}$ be a convex neighborhood of $x_0$. It is enough to prove that there exist $R > 0$ and $C > 0$ such that the norm of
  \[
  L_{y,x}\circ D\phi_{y}(z) - D\phi_{x}(z) : T_zM \ra T_xM,
  \]
  is bounded by $Cd(x, y)$ for every $x, y, z \in B\left(x_0, R\right)$. The function
  \[
  J : \mathcal{U} \times \mathcal{U} \times T\mathcal{U} \ra T\mathcal{U}
  \]
  defined by
  \[
  J\left(x, y,\left(z, h_{z}\right)\right) = L_{y,x}\circ D\phi_{y}(z)\left(h_z\right) - D\phi_{x}(z)\left(h_z\right)
  \]
  is $C^{\infty}$. Indeed, the map $\exp^{-1}$ is differentiable as a function defined on $\mathcal{U}\times \mathcal{U}$,
  and the parallel transport is the solution of an ordinary differential equation which
  depends $C^{\infty}$-wise on the initial data $x, y$, and consequently itself depends $C^{\infty}$-wise
  on $x, y$.

  Hence there exist $R > 0$ and $C_0 > 0$ such that $\norm DJ\left(x, y,\left(z, h_{z}\right)\right)\norm \leq C_0$, provided that $x, y, z \in B\left(x_0, R\right)$ and $\norm h_z\norm \leq R$, this implies that $\norm DJ\left(x, y,\left(z, h_{z}\right)\right)\norm \leq \frac{C_0}{R}=C$ provided that $x, y, z \in  B\left(x_0, R\right)$ and $\norm h_z\norm \leq 1$, and consequently
  \[
  \hspace{-3cm} d^{\ast}\left(J\left(x_1, y_1,\left(z_1, h_1\right)\right), J\left(x_2, y_2,\left(z_2, h_2\right)\right)\right)\leq
  \]
  \[
  \hspace{+2cm}C\left(d\left(x_1,x_2\right)+d\left(y_1,y_2\right)+d^{\ast}\left(\left(z_1,h_1\right), \left(z_2,h_2\right)\right)\right),
  \]
  provided that $x_i,y_i,z_i\in B\left(x_0, R\right)$ and $\norm h_i\norm\leq 1$, $ i = 1, 2$, where $d^{\ast}$
  is the distance induced in $TM$. Particularizing $x_1=x_2=y_2 = x$, $y_1 = y$, $z_1 = z_2 = z$ and
  $h_1=h_2=h$, we deduce
  \[
   \norm L_{y,x}\circ D\phi_{y}(z)\left(h\right) - D\phi_{x}(z)\left(h\right)\norm=\norm J\left(x, y,\left(z, h\right)\right) \norm \leq Cd(x,y),
   \]
  for every $x, y, z \in B\left(x_0, R\right)$ and $\norm h \norm \leq 1$. Thus
  \[
  \norm L_{y,x}\circ D\phi_{y}(z) - D\phi_{x}(z)\norm \leq Cd(x,y).
  \]

\end{proof}
%%%%%%%%%%%%%%%%%%%%%%%%%%%%%%%%%%%%%%%%%%%%%%%%%%%%%%%%%%

 We denote $\tilde{x} = P_S(x)$ and $\tilde{y} = P_S(y)$ for $x,y\in U$. We also shortened $L = L_{\tilde{y},\tilde{x}}$.
\begin{theorem}\label{Lemma2}
  If $S \subset M$ is $\varphi$-convex and $x_0\in S$, then there exist $R > 0$ and $C > 0$
  such that for every $0 < r < R$ and $x, y \in  B\left(x_0, r\right)$
  \begin{equation}\label{eq1}
    d\left(\tilde{x},\tilde{y}\right)\leq \frac{k\left(\tilde{x}\right)}{1-\beta-\alpha}d(x,y),
  \end{equation}
  where $k\left(\tilde{x}\right)$ is the Lipschitz constant of $\phi_{\tilde{x}}$ on a convex ball $B\left(\tilde{x}, \sigma\right)$ containing $x,y$, $\beta=\beta(x,y):=\varphi(\tilde{x})d\left(x,\tilde{x}\right)+\varphi(\tilde{y})d\left(y,\tilde{y}\right)$, and $\alpha=\alpha(y):=Cd\left(y,\tilde{y}\right)$.

\end{theorem}
\begin{proof}
 Let $x_0\in S$ and $R_0>0$, $C>0$ be the constants obtained from Lemma \ref{Lemma1} and let $R_0$ be small enough such that $B\left(x_0, R_0\right)$ is convex with compact closure and $B\left(x_0, R_0\right)\subset U$. We put $\bar{\varphi}:=\max_{x\in S\cap \overline{B\left(x_0, R_0\right)}}\varphi(x)$ and $\bar{r}:=\min_{x\in S\cap \overline{B\left(x_0, R_0\right)}}r(x)$ and then we define
  \[
  R:=\min\left\{\frac{R_0}{2},\frac{\bar{r}}{3}, \frac{1}{2\bar{\varphi}+C}\right\}.
  \]

  We now assume that $0 < r < R$ and $x, y \in  B\left(x_0, r\right)$. Hence $x,y,\tilde{x}, \tilde{y}\in B\left(x_0, R_0\right)$ and we have that
  \[
  \left\li \phi_{\tilde{x}}(x),\phi_{\tilde{x}}(\tilde{y})\right\ri_{\tilde{x}}\leq \varphi\left(\tilde{x}\right)\norm \phi_{\tilde{x}}(x)\norm d^2\left(\tilde{x},\tilde{y}\right),
  \]
  and
  \[
  \left\li \phi_{\tilde{y}}(y),\phi_{\tilde{y}}(\tilde{x})\right\ri_{\tilde{y}}\leq \varphi\left(\tilde{y}\right)\norm \phi_{\tilde{y}}(y)\norm d^2\left(\tilde{x},\tilde{y}\right).
  \]
  As inner product is invariant throughout parallel transport, we have also
  \[
  \left\li \phi_{\tilde{y}}(y),\phi_{\tilde{y}}(\tilde{x})\right\ri_{\tilde{y}}=\left\li L\left(\phi_{\tilde{y}}(y)\right),L\left(\phi_{\tilde{y}}(\tilde{x})\right)\right\ri_{\tilde{x}}.
  \]
  As $L\left(\phi_{\tilde{y}}(\tilde{x})\right) = - \phi_{\tilde{x}}(\tilde{y})$, we deduce
  \[
  \left\li L\left(\phi_{\tilde{y}}(y)\right),- \phi_{\tilde{x}}(\tilde{y})\right\ri_{\tilde{x}}\leq \varphi\left(\tilde{y}\right)\norm \phi_{\tilde{y}}(y)\norm d^2\left(\tilde{x},\tilde{y}\right),
  \]
  and consequently
  \begin{align*}
  \left\li \phi_{\tilde{x}}(x)-L\left(\phi_{\tilde{y}}(y)\right) ,\phi_{\tilde{x}}(\tilde{y})\right\ri_{\tilde{x}}&\leq \left(\varphi\left(\tilde{x}\right)\norm \phi_{\tilde{x}}(x)\norm+\varphi\left(\tilde{y}\right)\norm \phi_{\tilde{y}}(y)\norm \right) d^2\left(\tilde{x},\tilde{y}\right)\\
  &=\beta d^2\left(\tilde{x},\tilde{y}\right),
  \end{align*}
  where $\beta=\beta(x,y)=\varphi\left(\tilde{x}\right)\norm \phi_{\tilde{x}}(x)\norm+\varphi\left(\tilde{y}\right)\norm \phi_{\tilde{y}}(y)\norm$. Hence
  \[
  \beta d^2\left(\tilde{x},\tilde{y}\right)\geq \left\li \phi_{\tilde{x}}(x)-L\left(\phi_{\tilde{y}}(y)\right)-\phi_{\tilde{x}}(\tilde{y}) ,\phi_{\tilde{x}}(\tilde{y})\right\ri_{\tilde{x}}+ d^2\left(\tilde{x},\tilde{y}\right),
  \]
  since $\norm \phi_{\tilde{x}}(\tilde{y})\norm = d\left(\tilde{x},\tilde{y}\right)$. This implies
  \begin{align*}
    (1-\beta) d^2\left(\tilde{x},\tilde{y}\right)\leq & \left\li -\phi_{\tilde{x}}(x)+L\left(\phi_{\tilde{y}}(y)\right)+\phi_{\tilde{x}}(\tilde{y}) ,\phi_{\tilde{x}}(\tilde{y})\right\ri_{\tilde{x}} \\
    \leq & \norm -\phi_{\tilde{x}}(x)+L\left(\phi_{\tilde{y}}(y)\right)+\phi_{\tilde{x}}(\tilde{y})\norm d\left(\tilde{x},\tilde{y}\right),
  \end{align*}
  therefore
  \begin{align*}
    (1-\beta) d\left(\tilde{x},\tilde{y}\right)\leq  & \norm -\phi_{\tilde{x}}(x)+L\left(\phi_{\tilde{y}}(y)\right)+\phi_{\tilde{x}}(\tilde{y})\norm  \\
    \leq & \norm -\phi_{\tilde{x}}(x)+\phi_{\tilde{x}}(y)\norm+\norm -\phi_{\tilde{x}}(y)+L\left(\phi_{\tilde{y}}(y)\right)+\phi_{\tilde{x}}(\tilde{y})\norm \\
    \leq & k\left(\tilde{x}\right)d(x,y)+\norm -\phi_{\tilde{x}}(y)+L\left(\phi_{\tilde{y}}(y)\right)+\phi_{\tilde{x}}(\tilde{y})\norm,
  \end{align*}
  where $k\left(\tilde{x}\right)$ is the Lipschitz constant of $\phi_{\tilde{x}}$  on $B\left(\tilde{x}, 3r\right)\subset B\left(\tilde{x}, \bar{r}\right)$. Following with the second term, we
  have
  \begin{align*}
  \norm -\phi_{\tilde{x}}(y)+L\left(\phi_{\tilde{y}}(y)\right)+\phi_{\tilde{x}}(\tilde{y})\norm= & \norm \left(L\circ \phi_{\tilde{y}}-\phi_{\tilde{x}}\right)(y)-\left(L\circ \phi_{\tilde{y}}-\phi_{\tilde{x}}\right)(\tilde{y})\norm\\
  \leq & Ad\left(y,\tilde{y}\right),
  \end{align*}
  and consequently
  \[
  (1-\beta) d\left(\tilde{x},\tilde{y}\right)\leq k\left(\tilde{x}\right)d(x,y)+Ad\left(y,\tilde{y}\right),
  \]
  where $A$ is the Lipschitz constant of
  \[
  L\circ \phi_{\tilde{y}}-\phi_{\tilde{x}}: B\left(x_0, R\right)\ra T_{\tilde{x}}M.
  \]
 Using Lemma \ref{Lemma1}, $A \leq Cd\left(\tilde{x},\tilde{y}\right)$ and we
deduce
\[
(1-\beta) d\left(\tilde{x},\tilde{y}\right)\leq k\left(\tilde{x}\right)d(x,y)+\alpha d\left(\tilde{x},\tilde{y}\right),
\]
where $\alpha=Cd\left(y,\tilde{y}\right)$. Our choice of $R$ and $r$ implies that $1-\beta-\alpha>0$ and then we get the result. Indeed, we have
\begin{align*}
\beta+\alpha=&\varphi(\tilde{x})d\left(x,\tilde{x}\right)+\varphi(\tilde{y})
d\left(y,\tilde{y}\right)+Cd\left(y,\tilde{y}\right)\\
\leq & \bar{\varphi}d\left(x,x_0\right)+\bar{\varphi}d\left(y,x_0\right)+Cd\left(y,x_0\right)\\
\leq & \left(2\bar{\varphi}+C\right)r<1,
\end{align*}
since $x_0\in S$.

\end{proof}
%%%%%%%%%%%%%%%%%%%%%%%%%%%%%%%%%%%%%%%%%%%%%%%%%%%%%%
\begin{corollary}
  Let $S \subset M$ be $\varphi$-convex and $x_0\in S$. Then for every $\varepsilon > 0$ there exists $R > 0$ such that
  \[
   P_S : B\left(x_0, R\right) \ra S
   \]
   has Lipschitz constant less or equal than $1 + \varepsilon$.
\end{corollary}
%%%%%%%%%%%%%%%%%%%%%%%%%%%%%%%%%%%%%%%%%%%%%%%%%%%%%%%

Let $\gamma\in H^1(I,M)$ be such that $\gamma(t)\in S$ for all $t\in I$  and $V$ be a vector field along $\gamma$ such that $V\in H^1\left(I,\gamma^{-1}TM\right)$. We define
 \[
 P_{\gamma(t)}V(t):=P_{T^B_S\left(\gamma(t)\right)}V(t),\qquad \forall\ t\in I,
 \]
 where $T^B_S\left(\gamma(t)\right)$ is the Bouligand tangent cone to $S$ at the point $\gamma(t)$. Using \cite[Theorem 4.2]{Minimizing}, we have $P_{\gamma}V\in H^1\left(I,\gamma^{-1}TM\right)$ and
 \[
 \norm P_{\gamma(t)}V(t)\norm =\lim_{s\ra 0^{+}}\frac{d\left(\gamma(t), P_S\left(\exp_{\gamma(t)}sV(t)\right)\right)}{s}.
 \]

 We now consider a variation $\Gamma$ of $\gamma$ defined as follows
 \[
 \Gamma_s(t):=\exp_{\gamma(t)}\left(sV(t)\right), \quad \forall t\in I, s\geq 0,
 \]
   and for  sufficiently small s, we define $\tilde{\Gamma}_s(t):= P_S\left(\Gamma_s(t)\right)$ for all $t\in I$. Hence $\Gamma_s, \tilde{\Gamma}_s\in H^1(I,M)$.

To proceed, we need to estimate the following statement from below
\[
\liminf_{s\ra 0^{+}}\;\frac{f\left(\Gamma_s\right)-f\left(\tilde{\Gamma}_s\right)}{s}.
\]
Let $\im(\gamma)$ denote the image of $\gamma$  and $\rho:=\max_{t\in I}\norm V(t)\norm$.

%%%%%%%%%%%%%%%%%%%%%%%%%%%%%%%%%%%%%%%%%%%%%%%%%%%%%%%%%%%%%%%%%%%
\begin{theorem}\label{thm1}
Let $\gamma\in H^1(I,M)$ be such that $\gamma(t)\in S$ for all $t\in I$  and $V$ be a vector field along $\gamma$ such that $V\in H^1\left(I,\gamma^{-1}TM\right)$. Then there exists a piecewise constant function $\tau:I\ra \mathbb{R}$ such that
  \begin{equation*}
  \hspace{-4cm}\liminf_{s\ra 0^{+}}\;\frac{\frac{1}{2}\int_0^1\norm \dot{\Gamma_s}(t)\norm^2dt-\frac{1}{2}\int_0^1\norm \dot{\tilde{\Gamma}}_s(t)\norm^2dt }{s}\geq
\end{equation*}
\[
\hspace{3cm}-\int_{0}^{1}\left(\varphi(\gamma)+\tau\right)\norm V-P_{\gamma}V\norm \norm \dot{\gamma}\norm^2 dt.
\]
\end{theorem}
\begin{proof}
We first show that there exists a piecewise constant function $\tau$ on $I$ such that for all sufficiently small $s$ and for almost all $t\in I$,
\begin{equation}\label{q1}
  \norm \dot{\tilde{\Gamma}}_s(t)\norm\leq \frac{K(s)}{1-\left(2\varphi
\left(\tilde{\Gamma}_s(t)\right)+\tau(t)\right)d\left(\tilde{\Gamma}_s(t),\Gamma_s(t)\right)}\norm \dot{\Gamma_s}(t)\norm,
\end{equation}
 where $K$ is a function with the property that  $K\ra 1$ as $s\ra 0^{+}$.

 Indeed, using Theorem \ref{Lemma2} and the compactness of $\im(\gamma)$, there exist finitely many points $x_1,\ldots, x_m\in \im(\gamma)$ and some constants $R_i>0$, $C_i>0$  such that the inequality \eqref{eq1} holds on $B_i:=B\left(x_i,R_i\right)$ for $i=1,\ldots, m$ and the balls $B_i$, $i=1,\ldots , m$ cover $\im(\gamma)$.

 Let $t_i$, $i=0,\ldots, m$ be such that $t_0=0$, $t_m=1$, and $\gamma\left(\left[t_{i-1}, t_i\right]\right)\subset B_i$ for  $i=1,\ldots , m$. Since $\Gamma_s\ra \gamma$ uniformly on $I$, we can find $s_1>0$ small enough such that $s_1\rho<\min\left\{\frac{\bar{r}}{2}, \frac{R_i}{3}, i=1,\ldots , m\right\}$ and
 \[
 \Gamma_s(t)\in B_i, \quad \forall\ t\in \left[t_{i-1}, t_i\right], \quad \forall s<s_1,
 \]
 where $\bar{r}:=\min \left\{r(x): x\in S\cap \left(\bigcup \overline{B\left(x_i,2R_i\right)}\right)\right\}$.

 Let $s<s_1$ and $t\in \left(t_{j-1}, t_j\right)$ for some $1 \leq j\leq m$ be such that  $\dot{\Gamma_s}(t)$ and $\dot{\tilde{\Gamma}}_s(t)$ exist. Then for $h\in \mathbb{R}$ with sufficiently small $|h|$ we have  $t+h\in \left(t_{j-1}, t_j\right)$ and $\Gamma_s(t+h)\in B\left(\tilde{\Gamma}_s(t),s\rho\right)$, and hence using Theorem \ref{Lemma2} we obtain that
 \begin{equation}\label{q2}
d\left(\tilde{\Gamma}_s(t+h),\tilde{\Gamma}_s(t)\right)\leq \Theta\: d\left(\Gamma_s(t+h),\Gamma_s(t)\right),
\end{equation}
 where
 \[
 \Theta:=\frac{k\left(\tilde{\Gamma}_s(t)\right)}{1-\beta\left(\Gamma_s(t), \Gamma_s(t+h)\right)-C_jd\left(\tilde{\Gamma}_s(t+h),\Gamma_s(t+h)\right)},
 \]
 and $k$ is the Lipschitz constant of $\phi_{\tilde{\Gamma}_s(t)}$ on $B\left(\tilde{\Gamma}_s(t),s\rho\right)$. We now define the piecewise constant function $\tau$ on $I$ as $\tau(t)=C_i$ for all $t\in \left[t_{i-1}, t_i\right)$, $i=1,\ldots ,m$ and $\tau(t_m)=C_m$. Moreover,  let the sectional curvature of any plane of $M$  on $\bigcup B\left(x_i,R_i\right)$ be bounded  by $\Delta>0$, i.e. $|\se |\leq \Delta$, then putting
 \[
 K(s):=\frac{2s\rho\sqrt{\Delta}}{\sin \left(2s\rho\sqrt{\Delta}\right)},
 \]
 we have $k\left(\tilde{\Gamma}_s(t)\right)\leq K(s)$. Hence taking the limit of  \eqref{q2} as $h\ra 0$, we get the inequality \eqref{q1}.

%%%%%%%%%%%%%%%%%%%%%%%%%%%%%%%%%%%
For simplicity, we denote $x:=\Gamma_s(t)$ and $\tilde{x}:=\tilde{\Gamma}_s(t)$.  Using \eqref{q1}, we have

\[
\hspace{-6cm}  \frac{1}{2s}\int_0^1\left(\norm \dot{\Gamma_s}(t)\norm^2-\norm \dot{\tilde{\Gamma}}_s(t)\norm^2\right)dt\geq
  \]
\[
   \frac{1}{2s}\int_0^1\left(\norm \dot{\Gamma_s}(t)\norm^2 -\left(\frac{K(s)}{1-\left(2\varphi \left(\tilde{x}\right)+\tau(t)\right)d\left(\tilde{x},x\right)}\right)^2\norm \dot{\Gamma_s}(t)\norm^2\right)dt
\]
\begin{equation}\label{en1}
=\frac{1}{2s}\int_0^1\frac{1-K^2(s)}{\left[1-\left(2\varphi \left(\tilde{x}\right)+\tau\right)d\left(\tilde{x},x\right)\right]^2}\norm \dot{\Gamma_s}(t)\norm^2dt\;+
\end{equation}
\begin{equation}\label{en2}
\frac{1}{2s}\int_0^1\frac{\left(2\varphi \left(\tilde{x}\right)+\tau\right)d\left(\tilde{x},x\right)\left[\left(2\varphi \left(\tilde{x}\right)+\tau\right)d\left(\tilde{x},x\right)-2\right]}
{\left[1-\left(2\varphi \left(\tilde{x}\right)+\tau\right)d\left(\tilde{x},x\right)\right]^2}\norm \dot{\Gamma_s}(t)\norm^2dt.
\end{equation}

Note that
$\left(2\varphi\left(\tilde{x}\right)+c\right)d\left(\tilde{x},x\right)\ra 0$ uniformly on $I$ and  $\norm \dot{\Gamma_s}\norm^2\ra \norm \dot{\gamma}\norm^2$ in $L^1$ as $s\ra 0^{+}$. Moreover, we have $x,\tilde{x}\in B\left(\gamma(t),2s\rho\right)$ and then
\[
\frac{d\left(\tilde{\Gamma}_s(t),\Gamma_s(t)\right)}{s}\geq \frac{1}{K(s)}
\frac{\norm sV(t)-\exp_{\gamma(t)}^{-1}P_S\left(\exp_{\gamma(t)}sV(t)\right)\norm}{s},
\]
and
\[
\frac{1}{K(s)}
\frac{\norm sV(t)-\exp_{\gamma(t)}^{-1}P_S\left(\exp_{\gamma(t)}sV(t)\right)\norm}{s}\ra \norm V(t)-P_{\gamma(t)}V(t)\norm,
\]
in $L^1$ as $s\ra 0^{+}$. Therefore
\[
\lim_{s\ra 0^{+}}\eqref{en2}\geq -\int_{0}^{1}\left(2\varphi\left(\gamma\right)+\tau\right)\norm V-P_{\gamma}V\norm\norm \dot{\gamma}\norm^2dt.
\]

On the other hand, in \eqref{en1}, we have
\[
\frac{1}{\left[1-\left(2\varphi \left(\tilde{x}\right)+\tau\right)d\left(\tilde{x},x\right)\right]^2}\norm \dot{\Gamma_s}\norm^2\ra \norm \dot{\gamma}\norm^2,
\]
in $L^1$ and
\[
\frac{1-K^2(s)}{s}\ra 0
\]
 as $s\ra 0^{+}$ and consequently $\lim_{s\ra 0^{+}}\eqref{en1}=0$.
\end{proof}
%%%%%%%%%%%%%%%%%%%%%%%%%%%%%%%%%%%%%%%%%%%%%%%%%%%%%%%%%%%%%%%%%%%%%%%%%

In the sequel, we need to compute the derivative $\dot{\Gamma_s}(t)=\frac{\partial }{\partial t}\Gamma(s,t)$, where it exists. Let $t\in I$ be such that both $ \dot{\gamma}(t)$ and $D_tV(t)$ exist. Note that for fixed  $t$, $\Gamma(.,t)$ is a geodesic and $\Gamma$ is a variation of $\Gamma(.,t)$ among geodesics.  Let $(x^i)$ be the normal coordinates on $M$ centered at $\gamma(t)$ and $\left(x^i,v^i\right)$ be the corresponding coordinates  on $TM$. Then using the smooth approximations of $\gamma$ and $V$,  we have
\[
\dot{\Gamma_s}(t)=\left(\dot{\gamma}^i(t)+s\left(D_tV\right)^i(t)\right)\partial_i|_{\Gamma(s,t)},
\]
for sufficiently small $s$, where $\partial_i:=\partial/\partial x^i$ and $\dot{\gamma}^i(t)$, $\left(D_tV\right)^i(t)$ are the components of $\dot{\gamma}(t)$ and $D_tV(t)$ in these coordinates, respectively. Hence according to \cite{Klingenberg}, we obtain that
\[
\norm \dot{\Gamma_s}(t)\norm^2=\norm \dot{\gamma}(t)+sD_tV(t)\norm^2+O\left(s^2\right),
\]
where $O\left(s^2\right)/s\ra 0$ as $s\ra 0$.

%%%%%%%%%%%%%%%%%%%%%%%%%%%%%%%%%%%%%%%%%%%%%%%%%%%%%%%%%%%%%%%
\begin{theorem}\label{thm2}
  Let $\gamma\in \mathcal{A}$ and $\xi\in D^-f(\gamma)\subset L^2\left(I,\gamma^{-1}TM\right)$. Then there exists a piecewise constant function $\tau=\tau\left(\gamma\right)$ on $I$ such that for every vector field $V$ along $\gamma$ with the properties that $V\in H^1\left(I, \gamma^{-1}TM\right)$ and $V(0)=V(1)=0$,
  \[
  \int_{0}^{1}\li \dot{\gamma}, D_tV\ri dt\geq \int_{0}^{1}\li \xi, P_{\gamma}V\ri dt-\int_{0}^{1}\left(2\varphi(\gamma)+\tau\right)\norm V-P_{\gamma}V\norm \norm \dot{\gamma}\norm^2 dt.
  \]
\end{theorem}
\begin{proof}
  Note that
  \begin{align*}
    \li\dot{\gamma}(t), D_tV(t)\ri= & \frac{1}{2s}\norm \dot{\gamma}(t)+sD_tV(t)\norm^2-\frac{1}{2s}\norm \dot{\gamma}(t)\norm^2-\frac{s}{2}\norm D_tV(t)\norm^2 \\
    = & \frac{1}{s}\left(\frac{1}{2}\norm \dot{\Gamma_s}(t)\norm^2-\frac{1}{2}\norm \dot{\gamma}(t)\norm^2\right)+\left(\frac{O\left(s^2\right)}{s}-\frac{s}{2}\norm D_tV(t)\norm^2\right),
  \end{align*}
and hence
\[
\li\dot{\gamma}(t), D_tV(t)\ri=\lim_{s\ra 0^{+}}\frac{1}{s}\left(\frac{1}{2}\norm \dot{\Gamma_s}(t)\norm^2-\frac{1}{2}\norm \dot{\gamma}(t)\norm^2\right).
\]
Therefore we have
\[
\hspace{-6cm}\int_{0}^{1}\li \dot{\gamma}, D_tV\ri dt  -\int_{0}^{1}\li \xi, P_{\gamma}V\ri dt\,=
\]
\[
  \lim_{s\ra 0^{+}}\frac{1}{s}\int_{0}^{1}\left(\frac{1}{2}\norm \dot{\Gamma_s}(t)\norm^2-\frac{1}{2}\norm \dot{\gamma}(t)\norm^2-\left\li \xi , \exp_{\gamma(t)}^{-1}P_S\left(\exp_{\gamma(t)}sV(t)\right)\right\ri\right)dt
\]
\[
\geq \liminf_{s\ra 0^{+}}\frac{1}{s}\int_{0}^{1}\left(\frac{1}{2}\norm \dot{\tilde{\Gamma}}_s(t)\norm^2-\frac{1}{2}\norm \dot{\gamma}(t)\norm^2-\left\li \xi , \exp_{\gamma(t)}^{-1}\tilde{\Gamma}_s(t)\right\ri\right)dt
\]
\[
+\liminf_{s\ra 0^{+}}\frac{1}{2s}\int_{0}^{1}\left(\norm \dot{\Gamma_s}(t)\norm^2-\norm \dot{\tilde{\Gamma}}_s(t)\norm^2\right)dt.
\]
On the other hand, we have
\[
\liminf_{s\ra 0^{+}}\frac{1}{s}\int_{0}^{1}\left(\frac{1}{2}\norm \dot{\tilde{\Gamma}}_s(t)\norm^2-\frac{1}{2}\norm \dot{\gamma}(t)\norm^2-\left\li \xi , \exp_{\gamma(t)}^{-1}\tilde{\Gamma}_s(t)\right\ri\right)dt\geq
\]
\[
\hspace{-2.5cm}\liminf_{s\ra 0^{+}}\frac{f\left(\tilde{\Gamma}_s\right)-f\left(\gamma\right)-\left\li \xi , \exp_{\gamma}^{-1}\tilde{\Gamma}_s\right\ri_{L^2}}{s}=
\]
\[
\hspace{1cm}\liminf_{s\ra 0^{+}}\frac{f\left(\tilde{\Gamma}_s\right)-f\left(\gamma\right)-\left\li \xi , \exp_{\gamma}^{-1}\tilde{\Gamma}_s\right\ri_{L^2}}{\norm \exp_{\gamma}^{-1}\tilde{\Gamma}\norm_{L^2}}\times\frac{\norm \exp_{\gamma}^{-1}\tilde{\Gamma}\norm_{L^2}}{s}\geq 0,
\]
because $\xi\in D^-f(\gamma)$ and the function $\frac{\norm \exp_{\gamma}^{-1}\tilde{\Gamma}\norm_{L^2}}{s}$ is bounded by $2\norm V\norm_{L^2}$. Then Theorem \ref{thm1} implies that
 \[
  \int_{0}^{1}\li \dot{\gamma}, D_tV\ri dt- \int_{0}^{1}\li \xi, P_{\gamma}V\ri dt\geq -\int_{0}^{1}\left(2\varphi(\gamma)+\tau\right)\norm V-P_{\gamma}V\norm \norm \dot{\gamma}\norm^2 dt.
  \]

 \end{proof}
%%%%%%%%%%%%%%%%%%%%%%%%%%%%%%%%%%%%%%%%%%%%%%%%%%%%%%%%%%%%%

Let $H^1_0\left(I,\gamma^{-1} TM\right)$ denote the space of all proper vector fields  $V\in H^1\left(I, \gamma^{-1}TM\right)$ with $V(0)=V(1)=0$.
For $\gamma \in H^1(I, M)$, we can find an open neighborhood $W$ of $\im \left(\gamma\right)$ with compact closure. Let  $\delta\leq\se\leq\Delta$ on $U$, $\Delta\geq 0$ and $\bar{r}=\min_{t\in I}r\left(\gamma(t)\right)$. In the following lemma, we want to compute the covariant derivative of the $H^1$-vector field $\exp_{\gamma}^{-1}\eta$ along $\gamma$ with respect to $\dot{\gamma},\dot{\eta}$ for $\eta \in H^1(I, M)$ sufficiently near $\gamma$.

\begin{lemma}\label{lem2}
Suppose that $\eta \in H^1(I, M)$. If  $d_{\infty}\left(\gamma, \eta\right)< \bar{r}$, then for almost all $t\in I$,
\begin{equation}\label{e}
  \hspace{-2cm}\left\li \dot{\gamma}, D_t \left(\exp_{\gamma}^{-1}\eta\right)\right\ri\leq \left\li \dot{\gamma},  L_{\eta,\gamma}\dot{\eta}-\frac{1}{2}c\left(\gamma,\eta\right)\dot{\gamma}\right\ri
\end{equation}
\[
\hspace{4cm} +\frac{1}{2}|R|_{\infty} d^2\left(\gamma, \eta\right)\norm \dot{\gamma}\norm \norm \dot{\eta}\norm,
\]
where $c\left(\gamma,\eta\right)(t):=2\sqrt{\Delta}\:d\left(\gamma(t),\eta(t)\right)
\cot\left(\sqrt{\Delta}\:d\left(\gamma(t),\eta(t)\right)\right)$ for all $t\in I$ and $R$ denotes the curvature tensor on $M$.
\end{lemma}
\begin{proof}
  We define $V(t):=\phi_{\gamma(t)}\eta(t)$ for all $t\in I$. Then $V\in H^1\left(I,\gamma^{-1}TM\right)$ and for almost all $t\in I$,
  \begin{equation}\label{e1}
    D_tV(t)=D\phi_{\gamma(t)}\left(\eta(t)\right)\left(\dot{\eta}(t)\right)+\nabla V_1\left(\gamma(t)\right)\left(\dot{\gamma}(t)\right),
  \end{equation}
  where $V_1$ is a vector field on $B\left(\gamma(t),\bar{r}\right)$ defined by $V_1(x):= \exp_{x}^{-1}\left(\eta(t)\right)$. Thus $\nabla V_1\left(\gamma(t)\right)=\hess \left(-\frac{1}{2}d^2_{\eta(t)}\right)\left(\gamma(t)\right)$ and so using \cite[Lemma 3]{Minimizing}, we deduce that
  \begin{equation}\label{e2}
    \hess \left(-\frac{1}{2}d^2_{\eta(t)}\right)\left(\gamma(t)\right)\left(\dot{\gamma}(t)\right)^2\leq -\frac{1}{2}c\left(\gamma,\eta\right)(t)\norm \dot{\gamma}(t)\norm^2,
  \end{equation}
  where
  \[
   c\left(\gamma,\eta\right)(t)=2\sqrt{\Delta}\:d\left(\gamma(t),\eta(t)\right)
\cot\left(\sqrt{\Delta}\:d\left(\gamma(t),\eta(t)\right)\right).
   \]

   On the other hand, using \cite[p. 110]{Hardering}, we have
   \[
     \hspace{-3cm}\left\li \dot{\gamma}(t), D\phi_{\gamma(t)}\left(\eta(t)\right)\left(\dot{\eta}(t)\right)- L_{\eta(t),\gamma(t)}\dot{\eta}(t)\right\ri
     \]
   \[
       \leq\norm  D\phi_{\gamma(t)}\left(\eta(t)\right)- L_{\eta(t),\gamma(t)}\norm\norm \dot{\gamma}(t)\norm\norm \dot{\eta}(t)\norm
       \]
       \[
     \hspace{-1cm}\leq\frac{1}{2}|R|_{\infty}d^2\left(\gamma(t), \eta(t)\right)\norm \dot{\gamma}(t)\norm \norm \dot{\eta}(t)\norm.
 \]
  %Since $|R|_{\infty}\leq \frac{7}{3}\Delta$ on $W$,
 Hence we derive that
\begin{equation}\label{e3}
  \left\li \dot{\gamma}(t), D\phi_{\gamma(t)}\left(\eta(t)\right)\left(\dot{\eta}(t)\right)\right\ri\leq \left\li \dot{\gamma}(t), L_{\eta(t),\gamma(t)}\dot{\eta}(t)\right\ri+
 \end{equation}
 \[
  \frac{1}{2}|R|_{\infty} d^2\left(\gamma(t), \eta(t)\right)\norm \dot{\gamma}(t)\norm \norm \dot{\eta}(t)\norm.
\]

Therefore \eqref{e} is obtained from \eqref{e1},\eqref{e2} and \eqref{e3}.
\end{proof}
%%%%%%%%%%%%%%%%%%%%%%%%%%%%%%%%%%%%%%%%%%%%%%%%%%%%%%%%%%%%%%%%%%%

\begin{lemma}\label{lem3}
If $\eta \in H^1(I, M)$ and  $d_{\infty}\left(\gamma, \eta\right)< \bar{r}$, then for almost all $t\in I$,
\[
\left\li \exp_{\gamma}^{-1}\eta , D_t\left(\exp_{\gamma}^{-1}\eta\right)\right\ri=\left\li \exp_{\gamma}^{-1}\eta,  L_{\eta,\gamma}\dot{\eta}-\dot{\gamma}\right\ri.
\]
\end{lemma}
\begin{proof}
  Let $t\in I$ be such that $\dot{\gamma}(t)$, $\dot{\eta}(t)$ exist, then we have
  \[
   \hspace{-2cm} 2\left\li \left(\exp_{\gamma}^{-1}\eta\right)(t) , D_t\left(\exp_{\gamma}^{-1}\eta\right)(t)\right\ri= \frac{d}{dt}\left\norm\left(\exp_{\gamma}^{-1}\eta\right)(t)\right\norm^2
  \]
  \[
   \hspace{-3cm} = \frac{d}{dt}d^2\left(\gamma(t),\eta(t)\right)
  \]
  \[
   \hspace{1.7cm} =  \left\li -2\exp_{\gamma(t)}^{-1}\eta(t),\dot{\gamma}(t)\right\ri+\left\li -2\exp_{\eta(t)}^{-1}\gamma(t),\dot{\eta}(t)\right\ri
  \]
  \[
    =  2\left\li \exp_{\gamma(t)}^{-1}\eta(t),  L_{\eta(t),\gamma(t)}\dot{\eta}(t)-\dot{\gamma}(t)\right\ri.
    \]

\end{proof}
%%%%%%%%%%%%%%%%%%%%%%%%%%%%%%%%%%%%%%%%%%%%%%%%%%%%%%%%%%%%%%%%%%%%%%%

\begin{theorem}\label{thm4}
  Let $\gamma \in \mathcal{A}\cap W^{2,2}(I,M)$ and $\xi\in L^2\left(I,\gamma^{-1}TM\right)$ be such that
  \[
  \xi+D_t\dot{\gamma}\in N^P_S\left(\gamma\right),\ a.\ e.
   \]
   Then for all $\eta\in \mathcal{A}$ with the property that $d_{\infty}\left(\gamma, \eta\right)< \bar{r}$,
  \[
  \frac{1}{2}\int_{0}^{1}\norm \dot{\eta}(t)\norm^2dt\geq  \frac{1}{2}\int_{0}^{1}\left(c\left(\gamma,\eta\right)-1\right)\norm \dot{\gamma}(t)\norm^2dt+\int_{0}^{1}\left\li\xi,\exp_{\gamma}^{-1}\eta\right\ri dt
  \]
  \[
  -d_{\infty}(\gamma ,\eta )\norm\exp_{\gamma}^{-1}\eta\norm_{L^{2}}
  \big( \bar{\varphi}\norm\xi+
  D_t\dot{\gamma}\norm_{L^{2}}+\frac{1}{2}|R|_{\infty}\norm \dot{\gamma}\norm_{L^{\infty}}\norm\dot{\eta}\norm_{L^{2}}\big),
  \]
  where $\bar{\varphi}:=\max_{t\in I}\varphi\left(\gamma(t)\right)$.
\end{theorem}
\begin{proof}
  Let $\eta\in \mathcal{A}$ and $d_{\infty}\left(\gamma, \eta\right)< \bar{r}$, then we have
  \[
  \norm L_{\eta,\gamma}\dot{\eta}-\dot{\gamma}\norm^2=\norm \dot{\eta}\norm^2+\norm \dot{\gamma}\norm^2-2\left\li \dot{\gamma},L_{\eta,\gamma}\dot{\eta} \right\ri,
  \]
  and so
  \[
  \frac{1}{2}\norm \dot{\eta}\norm^2-\frac{c\left(\gamma,\eta\right)-1}{2}\norm \dot{\gamma}\norm^2=\frac{1}{2}\norm L_{\eta,\gamma}\dot{\eta}-\dot{\gamma}\norm^2+\left\li \dot{\gamma},L_{\eta,\gamma}\dot{\eta}-\frac{c\left(\gamma,\eta\right)}{2}\dot{\gamma} \right\ri.
  \]
  Moreover, using Lemma \ref{lem2} we have
  \[
  \frac{d}{dt}\left\li \dot{\gamma}, \exp_{\gamma}^{-1}\eta\right\ri \leq \left\li D_t\dot{\gamma}, \exp_{\gamma}^{-1}\eta\right\ri+\left\li \dot{\gamma}, L_{\eta,\gamma}\dot{\eta}-\frac{1}{2}c\left(\gamma,\eta\right)\dot{\gamma}\right\ri
  \]
  \[
  + \frac{1}{2}|R|_{\infty} d^2\left(\gamma, \eta\right)\norm \dot{\gamma}\norm \norm \dot{\eta}\norm.
  \]
  We now obtain that
  \[
  \frac{1}{2}\norm \dot{\eta}\norm^2-\frac{c\left(\gamma,\eta\right)-1}{2}\norm \dot{\gamma}\norm^2\geq \frac{1}{2}\norm L_{\eta,\gamma}\dot{\eta}-\dot{\gamma}\norm^2+\frac{d}{dt}\left\li \dot{\gamma}, \exp_{\gamma}^{-1}\eta\right\ri
  \]
  \[
  -  \left\li D_t\dot{\gamma}, \exp_{\gamma}^{-1}\eta\right\ri-\frac{1}{2}|R|_{\infty} d^2\left(\gamma, \eta\right)\norm \dot{\gamma}\norm \norm \dot{\eta}\norm.
  \]
  Hence by integrating from both side and noting that $\xi+D_t\dot{\gamma}\in N^P_S\left(\gamma\right),\ a.\ e.$, we have
  \[
  \frac{1}{2}\int_{0}^{1}\norm \dot{\eta}(t)\norm^2dt -  \frac{1}{2}\int_{0}^{1}\left(c\left(\gamma,\eta\right)-1\right)\norm \dot{\gamma}(t)\norm^2dt-\int_{0}^{1}\left\li\xi,\exp_{\gamma}^{-1}\eta\right\ri dt\geq
  \]
  \[
  \frac{1}{2}\int_{0}^{1}\norm L_{\eta,\gamma}\dot{\eta}-\dot{\gamma}\norm^2dt-\int_{0}^{1}\left\li\xi+D_t\dot{\gamma},\exp_{\gamma}^{-1}\eta\right\ri dt-\frac{1}{2}|R|_{\infty}\int_{0}^{1} d^2\left(\gamma, \eta\right)\norm \dot{\gamma}\norm \norm \dot{\eta}\norm dt
  \]
  \[
  \geq \frac{1}{2}\int_{0}^{1}\norm L_{\eta,\gamma}\dot{\eta}-\dot{\gamma}\norm^2dt-\int_{0}^{1}\varphi\left(\gamma\right)\norm\xi+
  D_t\dot{\gamma}\norm\norm\exp_{\gamma}^{-1}\eta\norm^2 dt
  \]
  \[
  -\frac{1}{2}|R|_{\infty}\int_{0}^{1} \norm \dot{\gamma}\norm \norm \dot{\eta}\norm\norm\exp_{\gamma}^{-1}\eta\norm^2 dt.
  \]
  On the other hand, using H\"{o}lder inequality we have
  \[
  \int_{0}^{1}\norm\xi+
  D_t\dot{\gamma}\norm\norm\exp_{\gamma}^{-1}\eta\norm^2 dt\leq d_{\infty}(\gamma ,\eta )
   \int_{0}^{1}\norm\xi+
  D_t\dot{\gamma}\norm\norm\exp_{\gamma}^{-1}\eta\norm dt
  \]
  \[
  \leq d_{\infty}(\gamma ,\eta)
  \norm\xi+ D_t\dot{\gamma}\norm_{L^{2}}\norm\exp_{\gamma}^{-1}\eta\norm_{L^{2}},
  \]
  and the following inequality,
  \begin{equation}\label{inft}
    \norm \dot{\gamma}\norm_{L^{\infty}}\leq \norm \dot{\gamma}\norm_{L^2}+\norm D_t\dot{\gamma}\norm_{L^1},
  \end{equation}
  implies that $\dot{\gamma}\in L^{\infty}(I,\gamma^{-1} TM)$. Therefore
  \[
  \int_{0}^{1} \norm \dot{\gamma}\norm \norm \dot{\eta}\norm\norm\exp_{\gamma}^{-1}\eta\norm^2 dt\leq \norm \dot{\gamma}\norm_{L^{\infty}}\norm\dot{\eta}\norm_{L^{2}}\norm\exp_{\gamma}^{-1}\eta\norm_{L^{2}}d_{\infty}(\eta ,\gamma),
  \]
  by H\"{o}lder inequality again, hence
  \[
  \frac{1}{2}\int_{0}^{1}\norm \dot{\eta}(t)\norm^2dt -  \frac{1}{2}\int_{0}^{1}\left(c\left(\gamma,\eta\right)-1\right)\norm \dot{\gamma}(t)\norm^2dt\geq
  \]
  \[
  \int_{0}^{1}\left\li\xi,\exp_{\gamma}^{-1}\eta\right\ri dt+\frac{1}{2}\int_{0}^{1}\norm L_{\eta,\gamma}\dot{\eta}-\dot{\gamma}\norm^2dt
  \]
  \[
  -d_{\infty}(\gamma ,\eta )\norm\exp_{\gamma}^{-1}\eta\norm_{L^{2}}\big(\bar{\varphi}\norm\xi+
  D_t\dot{\gamma}\norm_{L^{2}}
  +\frac{1}{2}|R|_{\infty}\norm \dot{\gamma}\norm_{L^{\infty}}\norm\dot{\eta}\norm_{L^{2}}\big).
  \]
  Then
  \[
  \frac{1}{2}\int_{0}^{1}\norm \dot{\eta}(t)\norm^2dt -  \frac{1}{2}\int_{0}^{1}\left(c\left(\gamma,\eta\right)-1\right)\norm \dot{\gamma}(t)\norm^2dt-\int_{0}^{1}\left\li\xi,\exp_{\gamma}^{-1}\eta\right\ri dt\geq
  \]
  \[
 -d_{\infty}(\gamma ,\eta )\norm\exp_{\gamma}^{-1}\eta\norm_{L^{2}}\big( \bar{\varphi}\norm\xi+
  D_t\dot{\gamma}\norm_{L^{2}}
  +\frac{1}{2}|R|_{\infty}\norm \dot{\gamma}\norm_{L^{\infty}}\norm\dot{\eta}\norm_{L^{2}}\big),
  \]
   and we get the result.
\end{proof}

%%%%%%%%%%%%%%%%%%%%%%%%%%%%%%%%%%%%%%%%%%%%%%%%%%%%%%%%%%%%%%%%%%%%%%%%%%%%%%%%%%%%%%%%%%%%%%%%%%%%%%%%%%%%%%%%%%%%
%%%%%%%%%%%%%%%%%%%%%%%%%%%%%%%%%%%%%%%%%%%%%%%%%%%%%%%%%%%%%%%%%%%%%%%%%%%%%%%%%%%%%%%%%%%%%%%%%%%%%%%%%%%%%%%%%%%%
\section{Weak geodesics as critical points of the energy functional}\label{sec4}
In this section, we characterize weak geodesics on $S$ as viscosity critical points of the energy functional.
\begin{theorem}\label{w1}
  Suppose that $\gamma\in \mathcal{A}$. Then $D^-f(\gamma)\neq \emptyset$ if and only if $\gamma\in W^{2,2}(I,M)$. Moreover, if $\xi\in L^2\left(I,\gamma^{-1}TM\right)$, then $\xi\in D^-f(\gamma)$ if and only if
  \[
  \xi(t)+D_t\dot{\gamma}(t)\in N^P_S\left(\gamma(t)\right),\ a.\ e. \ t\in I.
  \]
\end{theorem}
\begin{proof}
  Let $\xi\in D^-f(\gamma)$, then using Theorem \ref{thm2}, there exists a piecewise constant function $\tau$ on $I$ such that
  \begin{equation}\label{j2}
  \int_{0}^{1}\li \dot{\gamma}, D_tV\ri dt\geq \int_{0}^{1}\li \xi, P_{\gamma}V\ri dt-\int_{0}^{1}\left(2\varphi(\gamma)+\tau\right)\norm V-P_{\gamma}V\norm \norm \dot{\gamma}\norm^2 dt,
  \end{equation}
  for every  $V\in H_0^1\left(I, \gamma^{-1}TM\right)$.

  Since $T^B_S\left(\gamma(t)\right)$ is a closed convex cone for all $t$, we have
  $\norm  V(t)-P_{\gamma(t)}V(t)\norm\leq \norm V(t)\norm$ and $\norm V(t)\norm\geq \norm P_{\gamma(t)}V(t)\norm$ for all $t\in I$. Then similar to the proof of \cite[Lemma 3.5]{Canino2},  we obtain that
  \begin{equation}\label{j21}
  \int_{0}^{1}\li \dot{\gamma}, D_tV\ri dt\geq -\norm \xi\norm_{L^2}\norm V\norm_{L^2} -\int_{0}^{1}\left(2\varphi(\gamma)+\tau\right)\norm V\norm \norm \dot{\gamma}\norm^2 dt,
  \end{equation}
  by H\"{o}lder inequality.
  Therefore we have
 \begin{equation}\label{j3}
  \left|\int_{0}^{1}\li \dot{\gamma}, D_tV\ri dt\;\right|\leq \left(\norm \xi \norm_{L^2}+\left(2\bar{\varphi}+C\right)\norm \dot{\gamma} \norm_{L^2}^2\right)\norm V \norm_{L^{\infty}},
  \end{equation}
  for every $V\in H^1_0\left(I,\gamma^{-1} TM\right)$, where $C=\max_{I}\tau$.
  Hence taking a suitable sequence of vector fields $V_n\in H^1_0\left(I,\gamma^{-1} TM\right)$ in \eqref{j3} and passing to the limit as $n\ra \infty$, we conclude that
  \begin{equation}\label{j4}
  \norm \dot{\gamma}\norm_{L^{\infty}}\leq \norm \dot{\gamma}\norm_{L^1}+\norm \xi\norm_{L^2}+\left(2\bar{\varphi}+C\right)\norm \dot{\gamma} \norm_{L^2}^2,
  \end{equation}
  and hence $\dot{\gamma}\in L^{\infty}\left(I, \gamma^{-1} TM\right)$. Indeed, for given $t_0\in I$, let $w\in T_{\gamma(t_0)}M$ be  such that $\norm w\norm=1$ and $\li \dot{\gamma}(t_0),w\ri=\norm \dot{\gamma}(t_0)\norm$. We now consider a sequence of functions $u_n\in C^{\infty}_0\left(\mathbb{R}\right)$ with the properties that
  \[
  0\leq u_n\leq 1,\quad -tu'_n(t)\leq 1 \ \ \forall t\in \mathbb{R}, \quad \supp u_n\subseteq \left[-\frac{1}{n},\frac{1}{n}\right] \quad \forall n\in \mathbb{N}.
  \]
  These functions can be obtained as follows.  Let $\psi$ be the smooth function defined by
  \[
  \psi(t):=\left\{
  \begin{array}{ll}
    \lambda\exp\left(\frac{-1}{1-t^2}\right) & \mid t\mid<1 \\
    0 & \mid t\mid\geq 1,
  \end{array}\right.
  \]
  where $\lambda$ is the constant such that $\int \psi(t)dt=1$ and we put $u_n(t):=\frac{e}{\lambda}\psi(nt)$ for all $t\in \mathbb{R}$.\\
  We now define a sequence of  vector fields $V_n$ along $\gamma$ as
  \[
  V_n(t):=u_n\left(t-t_0\right)\left(t-t_0\right)L_{t_0,t}(w) \quad \forall\ t\in I,
  \]
  where $L_{t_0,t}$ denotes the parallel transport along $\gamma$ from $\gamma(t_0)$ to $\gamma(t)$. Then for $n$ sufficiently large, $V_n\in H^1_0\left(I,\gamma^{-1} TM\right)$ and $\norm V_n\norm_{L^{\infty}}\leq 1$. Putting $V_n$ in \eqref{j3}, we have
  \[
  \left\li \int_{0}^{1}u_n\left(t-t_0\right)L_{t,t_0}\dot{\gamma}(t)dt,w\right\ri+ \int_{0}^{1}u'_n\left(t-t_0\right)\left(t-t_0\right)\left\li L_{t,t_0}\dot{\gamma}(t),w\right\ri dt
  \]
  \[
  +\int_{0}^{1}u_n\left(t-t_0\right)\left(t-t_0\right)\left\li \dot{\gamma}(t), D_t\left(L_{t_0,t}(w)\right)\right\ri dt\leq \mathcal{C} \quad \forall n,
  \]
  where $\mathcal{C}:=\left(\norm \xi \norm_{L^2}+\left(2\bar{\varphi}+C\right)\norm \dot{\gamma} \norm_{L^2}^2\right)$. Hence passing to the limit as $n\ra \infty$, we conclude that $\norm \dot{\gamma}(t_0)\norm\leq \mathcal{C}+\norm \dot{\gamma}\norm_{L^1}$.

  On the other hand, using \eqref{j4} and by H\"{o}lder inequality, we have
  \begin{align*}
    \int_{0}^{1}\left(2\varphi(\gamma)+\tau\right)\norm V\norm \norm \dot{\gamma}\norm^2 dt \leq & \left(2\bar{\varphi}+C\right)\norm V \norm_{L^2}\left(\int_{0}^{1}\norm \dot{\gamma} \norm^4dt\right)^{1/2} \\
    \leq & \left(2\bar{\varphi}+C\right)\norm V \norm_{L^2}\norm\dot{\gamma} \norm_{L^2}\norm\dot{\gamma} \norm_{L^{\infty}} \\
    \leq & \left(2\bar{\varphi}+C\right)\norm V \norm_{L^2}\norm\dot{\gamma} \norm_{L^2}\\
     & \times \left(\norm \dot{\gamma}\norm_{L^1}+\norm \xi\norm_{L^2}+\left(2\bar{\varphi}+C\right)\norm \dot{\gamma} \norm_{L^2}^2\right).
  \end{align*}
  Hence applying this to \eqref{j21},  we get
   \begin{equation}\label{j5}
  \left|\int_{0}^{1}\li D_t\dot{\gamma}, V\ri dt\;\right|\leq \left(1+\left(2\bar{\varphi}+C\right)\norm \dot{\gamma}\norm_{L^2}\right)\left(\norm \xi\norm_{L^2}+\left(2\bar{\varphi}+C\right)\norm \dot{\gamma}\norm_{L^2}^2\right)\norm V \norm_{L^2},
  \end{equation}
  for every $V\in H^1_0\left(I,\gamma^{-1} TM\right)$.
  Hence taking a suitable sequence of vector fields $V_n\in H^1_0\left(I,\gamma^{-1} TM\right)$ in \eqref{j5} and passing to the limit as $n\ra \infty$, we obtain that
  \[
  \norm D_t\dot{\gamma}\norm_{L^2}\leq \left(1+\left(2\bar{\varphi}+C\right)\norm \dot{\gamma}\norm_{L^2}\right)\left(\norm \xi\norm_{L^2}+\left(2\bar{\varphi}+C\right)\norm \dot{\gamma}\norm_{L^2}^2\right).
  \]
  Indeed, for given $t_0\in I$, let $w\in T_{\gamma(t_0)}M$ be the vector such that $\norm w\norm=1$ and $\li D_t\dot{\gamma}(t_0),w\ri=\norm D_t\dot{\gamma}(t_0)\norm$. We now consider a sequence of functions $u_n\in C^{\infty}_0\left(\mathbb{R}\right)$ with the properties that
  \[
  0\leq u_n\leq 1,\quad \supp u_n\subseteq \left[-\frac{1}{n},\frac{1}{n}\right], \quad \forall n\in \mathbb{N},
  \]
  and we define a sequence of  vector fields $V_n$ along $\gamma$ as
  \[
  V_n(t):=u_n\left(t-t_0\right)L_{t_0,t}(w) \quad \forall\ t\in I.
  \]
   Then for $n$ sufficiently large, $V_n\in H^1_0\left(I,\gamma^{-1} TM\right)$ and $\norm V_n\norm_{L^2}\leq 1$. Putting $V_n$ in \eqref{j5}, we have
  \[
  \left\li \int_{0}^{1}u_n\left(t-t_0\right)L_{t,t_0}D_t\dot{\gamma}(t)dt,w\right\ri\leq \mathcal{D} \quad \forall n,
  \]
  where $\mathcal{D}:=\left(1+\left(2\bar{\varphi}+C\right)\norm \dot{\gamma}\norm_{L^2}\right)\left(\norm \xi\norm_{L^2}+\left(2\bar{\varphi}+C\right)\norm \dot{\gamma}\norm_{L^2}^2\right)$. Hence passing to the limit as $n\ra \infty$, we conclude that $\norm D_t\dot{\gamma}(t_0)\norm\leq \mathcal{D}$. Thus $D_t\dot{\gamma}\in L^2\left(I, \gamma^{-1} TM\right)$ and consequently  $\gamma\in W^{2,2}(I,M)$.

  We now show that
 \[
  D_t\dot{\gamma}(t)+\xi(t)\in N^P_S\left(\gamma(t)\right),\quad  \ a.\ e.\ t\in I.
  \]
  Using \cite[Lemma 1]{Minimizing}, we have $N^P_S(x)=\left(T^B_S(x)\right)^{\circ}$ and it suffices to prove that
  \[
  \left\li D_t\dot{\gamma}(t)+\xi(t),\omega\right\ri \leq 0 \qquad \forall\ \omega\in T^B_S\left(\gamma(t)\right), \quad  \ a.\ e.\ t\in I.
  \]

  Let $t_0\in I$ and $\omega\in T^B_S\left(\gamma(t_0)\right)$. Using a sequence of functions $u_n\in C^{\infty}_0\left(\mathbb{R}\right)$ with the properties that
  \[
  u_n\geq 0,\quad \supp u_n\subseteq \left[-\frac{1}{n},\frac{1}{n}\right], \quad \int u_n=1,\quad \forall n\in \mathbb{N}
  \]
  we construct a sequence of proper vector fields $V_n\in H^1_0\left(I,\gamma^{-1} TM\right)$ along $\gamma$ defined by
  \[
  V_n(t):=u_n\left(t-t_0\right)L_{t_0,t}(\omega) \quad \forall\ t\in I,
  \]
  where $L_{t_0,t}$ denotes the parallel transport along $\gamma$ from $\gamma(t_0)$ to $\gamma(t)$. Since $\gamma\in W^{2,2}(I,M)$, for very $n\in \mathbb{N}$ we have
  \begin{align*}
   \int_{0}^{1}\li \dot{\gamma}, D_tV_n\ri dt= & -\int_{0}^{1}\li D_t\dot{\gamma}, V_n\ri dt \\
    = & -\int_{0}^{1}\li u_n\left(t-t_0\right)L_{t,t_0}\left(D_t\dot{\gamma}\right), \omega\ri dt.
  \end{align*}
  Then from \eqref{j2} we obtain that
  \[
  -\left\li\int_{0}^{1} u_n\left(t-t_0\right)L_{t,t_0}\left(D_t\dot{\gamma}\right)dt, \omega\right\ri\geq \int_{0}^{1}u_n\left(t-t_0\right)\li \xi, P_{\gamma}\left(L_{t_0,t}\:\omega\right)\ri dt
  \]
  \[
  -\left(2\bar{\varphi}+C\right)\norm \dot{\gamma}\norm^2_{\infty}\int_{0}^{1}u_n\left(t-t_0\right)\norm L_{t_0,t}\:\omega-P_{\gamma}\left(L_{t_0,t}\:\omega\right)\norm  dt.
  \]
  Therefore when $n\ra \infty$, we derive that
  \[
  -\li D_t\dot{\gamma}(t_0),\omega\ri \geq \li\xi(t_0),P_{\gamma(t_0)}\:\omega\ri-\left(2\bar{\varphi}+C\right)\norm \dot{\gamma}\norm^2_{\infty}\norm \omega-P_{\gamma(t_0)}\:\omega\norm,
  \]
  and since $\omega\in T^B_S\left(\gamma(t_0)\right)$, $P_{\gamma(t_0)}\:\omega=\omega$ and then we get the result.

  For the converse, we assume that $\gamma\in W^{2,2}(I,M)$ and $\xi\in L^2\left(I,\gamma^{-1}TM\right)$ are such that
  \[
  \xi+D_t\dot{\gamma}\in N^P_S\left(\gamma\right),\ a.\ e.
  \]
  Then using Theorem \ref{thm4}, for every $\eta\in \mathcal{A}$ with $d_{\infty}\left(\gamma, \eta\right)< \bar{r}$ we have
  \[
  \frac{1}{2}\int_{0}^{1}\norm \dot{\eta}(t)\norm^2dt\geq  \frac{1}{2}\int_{0}^{1}\left(c\left(\gamma,\eta\right)-1\right)\norm \dot{\gamma}(t)\norm^2dt+\int_{0}^{1}\left\li\xi,\exp_{\gamma}^{-1}\eta\right\ri dt
  \]
  \[
   -d_{\infty}(\gamma ,\eta )\norm\exp_{\gamma}^{-1}\eta\norm_{L^{2}}
  \big( \bar{\varphi}\norm\xi+
  D_t\dot{\gamma}\norm_{L^{2}}+\frac{1}{2}|R|_{\infty}\norm \dot{\gamma}\norm_{L^{\infty}}\norm\dot{\eta}\norm_{L^{2}}\big),
  \]
  where $c\left(\gamma,\eta\right)=2\sqrt{\Delta}\:d\left(\gamma,\eta\right)
  \cot\left(\sqrt{\Delta}\:d\left(\gamma,\eta\right)\right)$.

  Since $c\left(\gamma,\eta\right)-2=O\big(d(\gamma ,\eta)^2\big)$, we have
  \[
   \Big|  \int_{0}^{1}\left(c\left(\gamma,\eta\right)-2\right)\norm \dot{\gamma}(t)\norm^2dt\Big| \leq
   K\int_0^1d(\gamma ,\eta )^2dt=K\left\norm \exp_{\gamma}^{-1}\eta\right\norm_{L^2}^2
   \]
  for a suitable constant $K$ that depends on the Taylor expansion of the tangent function at $0$ and $\norm \dot{\gamma}\norm_{L^{\infty}}^2$, and consequently,
  \[
    \liminf_{d_{\infty}\left( \eta,\gamma\right)\ra 0}\frac{f\left(\eta\right)-f\left(\gamma\right)-\left\li\xi,\exp_{\gamma}^{-1}\eta\right\ri_{L^2} }{\left\norm \exp_{\gamma}^{-1}\eta\right\norm_{L^2}}\geq 0.
   \]
   It follows that $\xi\in D^-f(\gamma)$.

   In particular, since
   \[
    -D_t\dot{\gamma}(t)+D_t\dot{\gamma}(t)=0\in N^P_S\left(\gamma(t)\right)\quad  \ a.\ e.\ t\in I,
   \]
   we deduce that $-D_t\dot{\gamma}\in D^-f(\gamma)$ and $D^-f(\gamma)\neq \emptyset$.
\end{proof}
%%%%%%%%%%%%%%%%%%%%%%%%%%%%%%%%%%%%%%%%%%%%%%%%%%%%%%%%%%%%%%%%%%

\begin{corollary}
  If $\gamma\in \mathcal{A}\cap W^{2,2}(I, M)$, then $-P_{\gamma}\left(D_t\dot{\gamma}\right)\in D^-f(\gamma)$ and
  \[
  \left\norm -P_{\gamma}\left(D_t\dot{\gamma}\right)\right\norm_{L^2} \leq \norm \xi \norm_{L^2}\quad \forall \, \xi\in D^-f(\gamma).
  \]
\end{corollary}
\begin{proof}
  Since $D_t\dot{\gamma}\in L^2\left(I,\gamma^{-1}TM\right)$, we also have $-P_{\gamma}\left(D_t\dot{\gamma}\right)\in L^2\left(I,\gamma^{-1}TM\right)$. Moreover,
  \[
  D_t\dot{\gamma}-P_{\gamma}\left(D_t\dot{\gamma}\right)\in N^P_{T^B_S\left(\gamma(t)\right)}(0)=N^P_S\left(\gamma(t)\right)\quad a.\ e. \ t\in I,
  \]
  because $T^B_S\left(\gamma(t)\right)$ is a closed convex subset of $T_{\gamma(t)}M$. Hence Theorem \ref{w1} implies that $-P_{\gamma}\left(D_t\dot{\gamma}\right)\in D^-f(\gamma)$.

  We now assume that $\xi\in D^-f(\gamma)$, thus $\xi+D_t\dot{\gamma}\in N^P_S\left(\gamma(t)\right)$ for almost all $t\in I$. It follows that
  \[
  \left\li \xi+D_t\dot{\gamma}, P_{\gamma}\left(D_t\dot{\gamma}\right)\right\ri\leq 0, \quad a.\ e.
  \]
  and hence
  \begin{align*}
    \left\li \xi, -P_{\gamma}\left(D_t\dot{\gamma}\right)\right\ri_{L^2}\geq & \left\li D_t\dot{\gamma}, P_{\gamma}\left(D_t\dot{\gamma}\right)\right\ri_{L^2} \\
    = & \left\li P_{\gamma}\left(D_t\dot{\gamma}\right), P_{\gamma}\left(D_t\dot{\gamma}\right)\right\ri_{L^2},
  \end{align*}
  that completes the proof.
\end{proof}

%%%%%%%%%%%%%%%%%%%%%%%%%%%%%%%%%%%%%%%%%%%%%%%%%%%%%%%%%%%%%%%%
\begin{theorem}
  Let $\gamma\in \mathcal{A}$, then $0\in D^-f(\gamma)$ if and only if $\gamma$ is a weak geodesic on $S$.
\end{theorem}
\begin{proof}
  Using Theorem \ref{w1}, we have $0\in D^-f(\gamma)$ if and only if $\gamma\in W^{2,2}(I,M)$ and $D_t\dot{\gamma}(t)\in N^P_S\left(\gamma(t)\right)$ for almost all $t\in I$.
\end{proof}
%%%%%%%%%%%%%%%%%%%%%%%%%%%%%%%%%%%%%%%%%%%%%%%%%%%%%%%%%%%%%%

\begin{proposition}
  If $\gamma\in \mathcal{A}$ is a weak geodesic on $S$, then $\gamma\in W^{2,\infty}\left(I, M\right)$ and $\gamma$ has constant speed.
\end{proposition}
\begin{proof}
   If $\gamma\in \mathcal{A}$ is a weak geodesic on $S$, then $0\in D^-f\left(\gamma\right)$ and hence using Theorem \ref{thm2}, there exists a piecewise constant function $\tau$ on $I$ such that for every vector field $V\in H^1_0\left(I, \gamma^{-1}TM\right)$,
  \[
  \int_{0}^{1}\li \dot{\gamma}, D_tV\ri dt\geq -\int_{0}^{1}\left(2\varphi(\gamma)+\tau\right)\norm V-P_{\gamma}V\norm \norm \dot{\gamma}\norm^2 dt.
  \]
  Then similar to the proof of Theorem \ref{w1} we deduce that $\dot{\gamma}\in L^{\infty}\left(I, \gamma^{-1}TM\right)$ and hence
  \[
  \left| \int_{0}^{1}\li \dot{\gamma}, D_tV\ri dt\right|\leq \left(2\bar{\varphi}+C\right)\norm \dot{\gamma}\norm^2_{L^{\infty}}\norm V\norm_{L^1},
  \]
  for all $V\in H^1_0\left(I, \gamma^{-1}TM\right)$. This implies that $D_t\dot{\gamma}\in L^{\infty}\left(I, \gamma^{-1}TM\right)$. Hence $\gamma\in W^{2,\infty}\left(I, M\right)$ and the function $\norm \dot{\gamma}\norm^2$ is Lipschitz on $I$. Indeed,  we have
  \[
   \left|\frac{d}{dt}\norm \dot{\gamma}(t)\norm^2\right|= \left|2\li \dot{\gamma}(t), D_t\dot{\gamma}(t)\ri\right|
     \leq  2\norm \dot{\gamma}\norm_{L^{\infty}}\norm D_t\dot{\gamma}\norm_{L^{\infty}}.
  \]
  Therefore similar to the proof of \cite[Theorem 3.8]{Canino2}, we show that $\frac{d}{dt}\norm \dot{\gamma}\norm^2=0$, a.e. on $I$. To this end, since $D_t\dot{\gamma}(t)\in N^P_S\left(\gamma(t)\right)$, a.e. on $I$, It suffices to prove that $\li \eta,\dot{\gamma}(t)\ri=0$ for all $\eta\in N^P_S\left(\gamma(t)\right)$.

  Let $\eta\in N^P_S\left(\gamma(t)\right)$ for some $t\in (0,1)$. Then for all $s>t$ and close enough to $t$ we have
  \[
  \left\li \eta,\frac{\exp_{\gamma(t)}^{-1}\gamma(s)}{s-t}\right\ri\leq \varphi\left(\gamma(t)\right)\norm \eta\norm \frac{d^2\left(\gamma(t),\gamma(s)\right)}{s-t}.
  \]
  Taking the limit as $s\ra t^{+}$, we conclude that $\left\li \eta,\dot{\gamma}(t)\right\ri\leq 0$. Similarly, we have $\left\li \eta,\dot{\gamma}(t)\right\ri\geq 0$ and then we get the result.

\end{proof}

%%%%%%%%%%%%%%%%%%%%%%%%%%%%%%%%%%%%%%%%%%%%%%%%%%%%%%%%%%%%%%%%%%%%%%%%%%%%%%%%%%%%%%%%%%%%%%%%%%%%%%%%%%%%%%%%%%%%%
%%%%%%%%%%%%%%%%%%%%%%%%%%%%%%%%%%%%%%%%%%%%%%%%%%%%%%%%%%%%%%%%%%%%%%%%%%%%%%%%%%%%%%%%%%%%%%%%%%%%%%%%%%%%%%%%%%%%%
%\begin{acknowledgements}
\noindent
\blu{\bf Acknowledgement}
The third-named author was supported by the Iran National Science Foundation
(INSF) under project No.4002602.
%\end{acknowledgements}

%\begin{acknowledgements}
%This work is based upon research funded by Iran National Science Foundation
%(INSF) under project No.4002602
%\end{acknowledgements}

%%%%%%%%%%%%%%%%%%%%%%%%%%%%%%%%%%%%%%%%%%%%%%%%%%%%%%%%%%%%%%%%%%%%%%%%%%%%%%%%%%%%%%%%%%%%%%%%%%%%%%%%%%%%

%%%%%%%%%%%%%%%%%%%%%%%%%%%%%%%%%%%%%%%%%%%%%%%%%%%%%%%%%%%%%%%%%%%%%%%%%%%%%%%%%%%%%%%%%%%%%%%%%%%%%%%%%%%%%%%%%%%%%%

\end{document}